\documentclass[11pt]{amsart}
\usepackage{amsmath,amssymb,amsthm,amscd}
\usepackage{amsmath}
\usepackage{amssymb}
\usepackage{amssymb,amsmath,amsthm,amsfonts,latexsym,euscript}
\usepackage{wasysym}
\usepackage{fancyhdr}
\fancypagestyle{plain}
\usepackage{hyperref}

\newtheorem{lemma}{Lemma}

\newtheorem*{example}{Example}
\newtheorem{theorem}{Theorem}

\usepackage{amssymb}
\usepackage[T1]{fontenc}
\usepackage[latin1]{inputenc}
\usepackage[usenames,dvipsnames]{xcolor}
\usepackage{graphicx}
\usepackage{setspace}
\usepackage{txfonts}
\usepackage{longtable}
\usepackage{blindtext}
\usepackage{ulem}
\usepackage[margin=1in]{geometry}

\AtEndDocument{{\footnotesize%
 \textsc{B.Rainear, Department of Mathematics, United States Naval Academy, Annapolis, MD., 21402} \par  
  \textit{E-mail address}, B.~Rainear: \texttt{benrainear@gmail.com} \par
  \textsc{K. Thompson, Department of Mathematics, United States Naval Academy, Annapolis, MD., 21402} \par
  \textit{E-mail address}, K. ~Thompson: \texttt{kthomps@usna.edu} \par
}}

\title{Elementary Proofs of Representation by Ternary Quadratic Forms}
\author{Benjamin Rainear}
\author{Katherine Thompson}
\begin{document}
\maketitle
\begin{abstract}
Mordell in 1958 \cite{Mordell} gave a new proof of the three squares theorem. Those techniques were generalized by Blackwell, et al., in 2016 \cite{BL} to characterize the integers represented by the remaining six ``Ramanujan-Dickson ternaries''. We continue the generalization of these techniques to four additional forms.
\end{abstract}
\section{Introduction}
The theory of quadratic forms has a long and rich history. Of particular interest is the question of representation of an integer by a form. In studying universal and almost universal positive definite forms (which, in particular, concerns four or more variables), knowing which integers are represented by ternary subforms not only is key from both theoretical and computational purposes, but also is a delicate and nontrivial matter. The four-squares theorem of Lagrange appeared in 1770, and the three-squares theorem of Legendre did not appear until 1797. And yet, assuming the three-squares theorem, the proof of the four squares theorem is at most a few lines. Even much more recent results such as the $451$ paper of Rouse \cite{Rouse}, which gives conditions under which quadratic forms are guaranteed to represent all odd positive integers, makes assumptions about certain ternary subforms--conditional on the Generalized Riemann Hypothesis; in considering the $24888$ escalators Rouse used knowledge of ternary subforms to handle $9812$ of these cases. Among these was a form where, if instead one took a more standard analytic approach and considered its corresponding theta series, would have required looking at a space of modular forms where the cuspidal subspace was $2604$ dimensional. \\
\\
This project is heavily influenced by recent work of the second author in Blackwell et al. \cite{BL}. Many of the results in that paper were not new; however, it was the technique that was unique. Concentrating mainly on the Ramanujan-Dickson ternaries (which, in particular, were ternary forms of determinant at most $10$), the authors showed which positive were represented by certain positive definite ternary forms. Proving what fails to be represented by a quadratic form is typically simple and straightforward; it is proving that $m \in \mathbb N$ is represented that is challenging. The authors began with a generic quadratic form of determinant $D$. They then showed a series of congruence conditions simultaneously held which guaranteed that the form in question represented a particular $m \in \mathbb N$, and also represented values that inequivalent forms of the same determinant $D$ failed to represent. \\
\\
As noted in the abstract, this paper begins by considering ternary forms not handled in \cite{BL}; indeed, the forms considered here are not as famous as the Ramanujan-Dickson forms and the following representation results do not seem to appear anywhere else in the literature. More crucially, we consider determinants $D$ much higher than those considered before, thus forcing many more candidate forms to be simultaneously eliminated. 
Last, we note that in \cite{BL} those excepted values by the other candidates happened to lie in the same congruence class, a congruence class which in turn had little to do directly with the determinant of the form. That is not the case here. The excepted values of the other candidates of determinant $D$ are of the form $Dk$, where $k$ is a quadratic nonresidue modulo $D$. While in \cite{BL} it was designed so that $x=1,y=z=0$ would be a vector so that $Q(x,y,z)$ produces an excepted value, that is not possible here because of what specifically is not represented by the other forms of determinant $D$. Therefore, the arguments for specific vector evaluation and subsequent elimination of other candidate forms are much more intricate. \\
\\
Our main results are:
\begin{theorem}
A positive integer $m$ is represented by $2x^2+2xy+2xz+2y^2+2yz+3z^2$ if and only if $m \neq 4^k (8 \ell +1)$.
\end{theorem}
\begin{theorem}
A positive integer $m$ is represented by $x^2+2y^2+2yz+ 6z^2$ if and only if $m \neq 4^k(8 \ell+5)$.
\end{theorem}
\begin{theorem}
A positive integer $m$ is represented by $x^2+3y^2+2yz+5z^2$ if and only if $m \neq 4^k(16 \ell+2)$.
\end{theorem}
\begin{theorem}
A positive integer $m$ is represented by $2x^2+2xy+3y^2+2yz+5z^2$ if and only if $m \neq 4^k(8 \ell+1)$.
\end{theorem}
The remainder of the paper is organized as follows: after a brief but more detailed background section, we proceed to the proofs of the theorems in order. As these proofs are constructive, we end with concrete examples.
\section{Background}
For general references on the theory of quadratic forms, we refer the reader to \cite{Cassels1} and \cite{Lam}.\\
\\
A \textbf{$n$-ary quadratic form} over $\mathbb Z$ is a polynomial $Q: \mathbb Z^n \to \mathbb Z$ given by $$Q(\vec{x}) = \displaystyle\sum_{1 \leq i \leq j \leq n}a_{ij}x_ix_j \in \mathbb Z[x_1,...,x_n].$$ We say a quadratic form is \textbf{positive definite} if $Q(\vec{x}) \geq 0$ for all $\vec{x} \in \mathbb Z^n$ and if $Q(\vec{x})=0$ if and only if $\vec{x} = \vec{0}$.\\
\\
To each quadratic form there is an associated symmetric matrix $A_Q\in \mathcal{M}_n \left( \tfrac{1}{2} \mathbb Z \right)$ whose entries $c_{ij}$ are given by $$c_{ij}= \begin{cases} a_{ij}, & i=j \\ a_{ij}/2, & i \neq j. \end{cases}$$ With this, we note that $\vec{x}\phantom{l}^T A_Q \vec{x} = Q(\vec{x})$. When $A_Q \in \mathcal{M}_n(\mathbb Z)$ we say that $Q$ is \textbf{classically integral}. We say that the \textbf{determinant} $\det(Q) = \det(A_Q)$. Last, we say that two forms $Q$ and $Q'$ are \textbf{equivalent} if there exists $M \in SL_n(\mathbb Z)$ such that $M^tA_{Q}M = A_{Q'}$.\\
\\
From now on, when we use the word ``form'' we mean ``ternary, positive definite, classically integral quadratic form.''\\
\\
In \cite{BL} the main idea was as follows: to show that some square-free $m \in \mathbb N$ is represented by a form $R$ of determinant $D$, suppose $$mR(x,y,z) = (Ax+By+mz)^2 + ax^2+2hxy+by^2.$$ If one can show all of the following conditions hold for some integers $A,B,h,b,a$:
\begin{itemize}
\item $ab-h^2=Dm$;
\item $A^2+a\equiv B^2+b \equiv 2AB+2h \equiv 0 \pmod{m}$;
\item $\left( \dfrac{-Dm}{a}\right) = 1 = \left( \dfrac{-a}{p} \right)$ where $p \vert m$ is prime;
\item $R(\vec{x})=k$ where $k \in \mathbb N$ is not represented by forms $Q'$ of determinant $D$ with $Q'$ not equivalent to $Q$,
\end{itemize}
then $R$ must be equivalent to $Q$. We continue to use this same basic approach. However, in \cite{BL} simplicities were made that can no longer be afforded. Namely, they took $b \equiv B \equiv h \equiv 0 \pmod{m}$, and they made $\tfrac{A^2+a}{m}=k$ (therefore, $\vec{x}=(1,0,0)$). Here, we do not (necessarily) have $b \equiv B \equiv h \equiv 0 \pmod{m}$, and moreover, we have (with one exception) $\vec{x} = (1,1,0)$. \\
\\
We note that the choice of $b \equiv B \equiv h \equiv 0 \pmod{m}$ in \cite{BL} was to make $2AB+2h \equiv 0 \pmod{m}$ immediate. One key realization here is that once $A$ and $B$ have been determined so that $A^2+a \equiv B^2+b \equiv 0 \pmod{m}$, even if $A,B \not \equiv 0 \pmod{m}$, one still can ensure $2AB+2h \equiv 0 \pmod{m}$ and in fact that $AB+h \equiv 0 \pmod{m}$. For a prime $p \vert m$: 
\begin{eqnarray*}
(AB+h)^2 & \equiv & A^2B^2 +2ABh+h^2 \pmod{p}\\
(AB+h)^2 & \equiv & h^2+2ABh+h^2 \pmod{p} \\
(AB+h)^2-2h(h+AB) & \equiv & 0 \pmod{p}\\
(AB+h)(AB-h) & \equiv & 0 \pmod{p}.
\end{eqnarray*}
And so, $AB+h \equiv 0 \pmod{p}$ can be ensured for all $p \vert m$, which by the Chinese Remainder Theorem guarantees $AB+h \equiv 0 \pmod{m}$ has a solution.
\section{Proof of Theorem 1}
In this section, we provide a proof of Theorem $1$, noting that there are three forms of determinant $7$: $Q_1: 2x^2+2xy+2xz+2y^2+2yz+3z^2$, $Q_2: x^2+y^2+7z^2$, $Q_3: x^2+2y^2+4z^2+2yz$.
\begin{lemma}
For any $m \equiv 1 \pmod{8}$, $Q_1$ does not represent $m$.
\end{lemma}
\begin{proof}
This is a simple exercise left to the reader.
\end{proof}
\begin{lemma}
If $m$ is odd, $4m$ is represented by $Q_1$ if and only if $m$ is.
\end{lemma}
\begin{proof}
One direction is trivial. So suppose $4m$ is represented, where $m$ is odd. Note that then $4m \equiv 4 \pmod{8}$. By a computer search, one can determine that this forces all of $x,y,z$ to be even. Substituting $x=2X, y=2Y, z=2Z$ we have $$4m = 4(2X^2)+4(2XY)+4(2XZ)+4(2Y^2)+4(2YZ)+4(3Z^2)$$ and dividing through, we see $m$ is represented by $Q_1$.
\end{proof}
\begin{lemma}
If $m = 4^k (8 \ell +1)$ for integers $k, \ell$, then $m$ is not represented by $Q_1$.
\end{lemma}
\begin{proof}
This follows immediately from the previous two lemmas. 
\end{proof}
\begin{lemma}
If $m=7n$ where $n \equiv 3,5,6 \pmod{7}$, then $m$ is not represented by $Q_2$.
\end{lemma}
\begin{proof}
Suppose that $Q_2$ represents $7n$ for some $n \in \mathbb Z$. Then necessarily $x \equiv y \equiv 0 \pmod{7}$, and substituting $x=7X$ and $y=7Y$ we see 
\begin{eqnarray*}
49X^2+49Y^2+7z^2 & = & 7n \\
7X^2+7Y^2+z^2 & = & n.
\end{eqnarray*}
This implies that $n$ is a quadratic residue modulo $7$.
\end{proof}
\begin{lemma}
If $m=7n$ where $n \equiv 3,5,6 \pmod{7}$ then $m$ is not represented by $Q_3$.
\end{lemma}
\begin{proof}
A computer search shows that for $Q_3$ to represent any number divisible by $7$, $x \equiv 0 \pmod{7}$. Moreover, modulo $7$, $(y,z) \in \{ \pm (1,5), \pm (2,3), \pm (3,1)\}$. Consider the first case, writing $x=7X$, $y=7Y+1$ and $z=7Z+5$:
\begin{eqnarray*}
	49X^2+2(7Y+1)^2+4(7Z+5)^2+2(7Y+1)(7Z+5) & = & 7n \\
	7X^2+14Y^2+14Y+14YZ+28Z^2+42Z+16 & = & n.
\end{eqnarray*}
This implies that $n \equiv 2 \pmod{7}$. Similarly, substituting $y=7Y+2, z = 7Z+3$ we have
\begin{eqnarray*}
49X^2+2(7Y+2)^2 + 4(7Z+3)^2+2(7Y+2)(7Z+3) & = & 7n \\
7X^2+14Y^2+14Y+14YZ+28Z^2+28Z+8 & = & n,
\end{eqnarray*}
which means $n \equiv 1 \pmod{7}$. Last, with $y=7Y+3, z = 7Z+1$ we have
\begin{eqnarray*}
49X^2+2(7Y+3)^2+4(7Z+1)^2+ 2(7Y+3)(7Z+1) & = & 7n \\
7X^2+14Y^2+14YZ+28Z^2+14Z+4 & = & n,
\end{eqnarray*}
and $n \equiv 4 \pmod{7}$.
\end{proof}
Now we suppose $m \not\equiv 1 \pmod{8}$ is squarefree. We will show $m$ is represented by $Q_1$. We proceed by cases. In the interest of space, and so as not to belabor the reader with repetition, we make note of when cases become identical to those completed in more detail.
\begin{itemize}

\item[(Case 1)] Suppose $m \equiv 3 \pmod{4}$. We choose $a \equiv 1 \pmod{4}$ and $a \equiv 3 \pmod{49}$ a prime such that $(\tfrac{-a}{p})=1$ for all primes $p \vert m$. This guarantees $$\left( \frac{-7m}{a} \right) = \left( \dfrac{-1}{a}\right) \left( \dfrac{7}{a}\right) \left( \dfrac{m}{a}\right) = 1 = \left( \dfrac{-a}{m} \right).$$
Considering the equation $ab-h^2=7m$, we see that modulo $7$, $(b,h) \in \{ (0,0), (3, \pm 3), (5, \pm 1), (6, \pm 2)\}$. Switching $h$ with $-h$ as necessary, we can safely assume modulo $7$, $(b,h) \in \{  (0,0), (3,4), (5,1), (6,2)\}$. This automatically guarantees 
there is a solution to $$(A+B)^2 + a+b+2h \equiv 0 \pmod{7}$$ which means that $Q$ will represent multiples of seven. Moreover, for all but the case where $(b,h) = (3,4)$ there are solutions to each of $$(A+B)^2 + a+b+2h \equiv 0,7,14,21,28, 35,42 \pmod{49}$$ This will guarantee that when $x \equiv y \equiv 1 \pmod{49}$ and $z =0$, $Q$ will represent $7n$ where $n$ is a quadratic nonresidue modulo $7$ (setting the equation to $21,35,42 \pmod{42}$ when $m$ is a quadratic residue $\pmod{7}$ and to $7,14,28$ when $m$ is a quadratic nonresidue $\pmod{7}$ ).  In the case where $(b,h) = (3,4)$ there is a solution to $$4A^2+4AB+B^2 +4a+b+4h \equiv 2^2(A^2+a) + 2(2AB+2h)+B^2+b \equiv 0,7,14, 21, 28, 35, 42 \pmod{49}$$ which means when $x \equiv 2 \pmod{49}, y \equiv 1 \pmod{49}$ and $z=0$ then $Q$ will represent $7n$ where $n$ is a quadratic nonresidue modulo $7$ (with similar restrictions based on $m$ being a quadratic residue $\pmod{7}$ or not.). Regardless, in each case $A$ and $B$ are predetermined $\pmod{m}$, such that $A^2+a \equiv B^2 + b \equiv 2(AB+h) \equiv 0 \pmod{m}$. 


\item[(Case 2)] Suppose $m \equiv 6 \pmod{8}$. Then $m=2m'$ where $m' \equiv 3 \pmod{4}$. Choose $a \equiv 1 \pmod{8}$ and $a \equiv 3 \pmod{49}$ to be prime, where additionally $(\tfrac{-a}{p})=1$ for all primes $p \vert m'$. This guarantees $$\left( \frac{-7m}{a} \right) = \left( \dfrac{-1}{a}\right) \left( \dfrac{7}{a}\right) \left( \dfrac{2}{a} \right) \left( \dfrac{m'}{a}\right) = 1= \left( \dfrac{-a}{m} \right).$$ The rest of this case is identical to (Case 1).


\item[(Case 3)] Suppose $m \equiv 5 \pmod{8}$. Let $a=2a'$ where $a'$ is a prime satisfying $a' \equiv 26 \pmod{49}$, $a' \equiv 1 \pmod{4}$ and where for all $p \vert m$, $(\tfrac{-2a'}{p})=1$. This guarantees $$\left( \dfrac{-7m}{a'} \right) = \left( \dfrac{-1}{a'}\right) \left( \dfrac{7}{a'}\right) \left( \dfrac{m}{a'}\right) =1 = \left( \dfrac{-a}{m}\right).$$ Moreover, if $a' \equiv 26 \pmod{49}$, $2a' \equiv 3 \pmod{49}$. That means the rest of this case will reduce to (Case 1).


\item[(Case 4)] Suppose $m \equiv 2 \pmod{8}$. Then $m=2m'$ where $m' \equiv 1, 5 \pmod{8}$. Let $a \equiv 5 \pmod{8}$ and $a \equiv 3 \pmod{49}$ be prime with $(\tfrac{-a}{p})=1$ for all primes $p \vert m'$. This guarantees  $$\left( \frac{-7m}{a} \right) = \left( \dfrac{-1}{a}\right) \left( \dfrac{7}{a}\right) \left( \dfrac{2}{a} \right) \left( \dfrac{m'}{a}\right) = 1.$$ The rest of this case is identical to (Case 1).

\end{itemize}

\section{Proof of Theorem 2}
There are three forms of determinant $11$: $Q_1: x^2+2y^2+2yz+6z^2$, $Q_2: x^2+y^2+11z^2$, and $Q_3: x^2+3y^2+2yz+4z^2$.
\begin{lemma}
If $m \equiv 5 \pmod{8}$, then $Q_1$ does not represent $m$.
\end{lemma}
\begin{proof}
This is a simple proof by exhaustion and is left to the reader.
\end{proof}
\begin{lemma}
A natural number $m \in \mathbb N$ is represented by $Q_1$ if and only if $4m$ is.
\end{lemma}
\begin{proof}
One direction is trivial. Suppose $4m$ is represented by $Q_1$. Looking $\pmod{2}$ this implies $x$ is even. Looking then $\pmod{4}$, we see that $y$ and $z$ cannot both be odd; however, any one of them even forces the other to be even. Writing $x=2X, y=2Y, z=2Z$ and dividing through by $4$ gives the result.
\end{proof}
\begin{lemma}
If $m=4^k(8 \ell +5)$ then $m$ is not represented by $Q_1$.
\end{lemma}
\begin{proof}
This follows immediately from the previous lemmas. 
\end{proof}
\begin{lemma}
If $m=11n$ where $n$ is a quadratic nonresidue modulo $11$, then $m $ is not represented by $Q_2$.
\end{lemma}
\begin{proof}
Without loss of generality, suppose $m$ is squarefree. Suppose $Q_2$ represents $m=11n$. This immediately forces $x \equiv y \equiv 0 \pmod{11}$. Substituting $x=11X$ and $y=11Y$ this means
\begin{eqnarray*}
121X^2+121Y^2+11z^2 & = & 11n\\
11X^2+11Y^2+z^2 & = & n
\end{eqnarray*}
which means that $n$ is a quadratic residue modulo $11$.
\end{proof}
\begin{lemma}
If $m=11n$ where $n$ is a quadratic nonresidue modulo $11$, then $m$ is not represented by $Q_3$.
\end{lemma}
\begin{proof}
For $Q_3$ to represent any multiple of $11$, we must have $x \equiv 0 \pmod{11}$. Assuming that $m$ is squarefree, this moreover, modulo $11$, means $(y,z) \in \{ \pm (1,8), \pm (2,5), \pm (3,2), \pm (4,10), \pm (5,7)\}$. Again, we proceed by cases. Writing first $y=11Y+1$ and $z=11Z+8$ gives
\begin{eqnarray*}
121X^2+3(11Y+1)^2+4(11Z+8)^2 + 2(11Y+1)(11Z+8) & = & 11n \\
11X^2 + 33Y^2+44Z^2+22Y + 66Z+22YZ+25 & = & n
\end{eqnarray*}
which means that $n \equiv 3 \pmod{11}$.\\
Next if $y=11Y+2$, $z=11Z+5$ we have
\begin{eqnarray*}
121X^2+3(11Y+2)^2+4(11Z+5)^2 + 2(11Y+2)(11Z+5) & = & 11n \\
11X^2+33Y^2+44Z^2+22Y+44Z+22YZ + 12 & = & n
\end{eqnarray*}
which means $n \equiv 1 \pmod{11}$. Supposing next $y=11Y+3$ and $z=11Z+2$ we see
\begin{eqnarray*}
121X^2+3(11Y+3)^2+4(11Z+2)^2 + 2(11Y+3)(11Z+2) & = & 11n \\
11X^2+33Y^2+44Z^2+22Y+22Z+22YZ + 5 & = & n
\end{eqnarray*}
and again $n$ is a quadratic residue modulo $11$. When $y=11Y+4$ and $z=11Z+10$ we have
\begin{eqnarray*}
121X^2+3(11Y+4)^2+4(11Z+10)^2 + 2(11Y+4)(11Z+10) & = & 11n \\
11X^2+33Y^2+44Z^2+44Y+88Z+22YZ + 48 & = & n
\end{eqnarray*}
which makes $n \equiv 4 \pmod{11}$. Last, if $y=11Y+5$ and $z = 11Z+7$:
\begin{eqnarray*}
121X^2+3(11Y+5)^2+4(11Z+7)^2 + 2(11Y+5)(11Z+7) & = & 11n \\
11X^2+33Y^2+44Z^2+44Y+66Z+22YZ + 31& = & n
\end{eqnarray*}
and $n \equiv 9 \pmod{11}$.
\end{proof}
Next, suppose that $m \neq 5 \pmod{8}$ is squarefree. This gives the following cases:
\begin{itemize}
\item[(Case 1)] $m \equiv 1 \pmod{8}$. We set $a = 2a_1$, where $a_1$ is an prime 
satisfying $a_1 \equiv 1 \pmod 4$ and $a_1 \equiv 1 \pmod{121}$ and  $(\frac{-a}{p}) = 1$,  for all primes $ p \vert m$. This guarantees $$\left( \dfrac{-11m}{a} \right) = \left( \dfrac{-a}{m} \right)=1.$$Replacing $h$ with $-h$ as necessary, we can assume $(b,h) \pmod{11} \in \{(0,0), (2,2), (6,1), (7,5), (8,4), (10,3)\}$. This automatically guarantees a solution to $$(A+B)^2=a+b+2h \equiv 0 \pmod{11}$$ which means that $Q$ will represent multiples of $11$. Moreover, for each pair $(b,h)$ there are solutions to $$(A+B)^2+a+b+2h \equiv 11k \pmod{121}$$ for $k=0,1,2,...,10$. This will guarantee that when $x \equiv y \equiv 1 \pmod{121}$ and $z=0$, $Q$ will represent $11k$ where $k$ is a quadratic nonresidue modulo $11$ (setting the equation to $11,33,44,55,99$ when $m$ is a quadratic residue $\pmod{11}$ and to $22,55,66,77,88,110$ when $m$ is a quadratic nonresidue $\pmod{11}$).


\item[(Case 2)] $m \equiv 2 \pmod{8}$. We start with writing $m = 2\ell$, where $\ell \equiv 1 \pmod{4}$. To ensure $$\left( \dfrac{-11m}{a}\right) =\left( \dfrac{-a}{m} \right)=1$$ we set $a$ to be a prime satisfying $a \equiv 5 \pmod{8}$, $a \equiv 2 \pmod{121}$, and $(\tfrac{-a}{p})=1$ for all odd $p \vert m$. The rest of the proof then follows (Case $1$).

\item[(Case 3)] $m \equiv 3 \pmod{4}$. Here we choose a prime $a \equiv 1 \pmod{4}$, $a \equiv 2 \pmod{121}$ and $(\tfrac{-a}{p})=1$ for all primes $p \vert m$. Then $$\left( \dfrac{-11m}{a}\right) = \left( \dfrac{-a}{m}\right) = 1$$ and the rest of the proof mimics (Case $1$).

\item[(Case 4)] $m \equiv 6 \pmod{8}$. We write $m = 2\ell$, where $\ell \equiv 3 \pmod 4$. Here we choose $a$ prime satisfying $a \equiv 1 \pmod{8}$, $a \equiv 2 \pmod{121}$, and $(\tfrac{-a}{p})=1$ for all primes $p \vert m$. Then $$\left( \dfrac{-11m}{a} \right) = \left( \dfrac{-a}{m} \right) = 1$$ and the rest of the proof follows like the others. 
%
%
\end{itemize}

\section{Proof of Theorem 3}
There are three forms of determinant $14$: $Q_1 : x^2+3y^2+2yz+5z^2$, $Q_2: x^2+y^2+14z^2$, $Q_3: x^2+2y^2+7z^2$.
\begin{lemma}
If $m \equiv 2 \pmod{16}$ then $m$ is not represented by $Q_1$.
\end{lemma}
\begin{proof}
We leave the proof to the reader. 
\end{proof}
\begin{lemma}
Let $m$ be even. Then $m$ is represented by $Q_1$ if and only if $4m$ is represented by $Q_1$.
\end{lemma}
\begin{proof}
For the nontrivial direction, if $m$ is even, then $4m \equiv 0 \pmod{8}$, which forces $x,y,z$ all even. 
\end{proof}
\begin{lemma}
If $m=4^k(16 \ell+2)$, then $m$ is not represented by $Q_1$.
\end{lemma}
\begin{proof}
This follows from the previous lemmas. 
\end{proof}
\begin{lemma}
If $m=7n$ where $n$ is a quadratic nonresidue modulo $7$, then $m$ is not represented by $Q_2$.
\end{lemma}
\begin{proof}
Suppose $m=7n$ is squarefree and is represented by $Q_2$. This forces $x \equiv y \equiv 0 \pmod{7}$. Substituting $x=7X$ and $y=7Y$ and simplifying gives $$7X^2+7Y^2+2z^2 = n$$ which means that $n$ is a quadratic residue modulo $7$.
\end{proof}
\begin{lemma}
If $m=7n$ where $n$ is a quadratic nonresidue modulo $7$, then $m$ is not represented by $Q_3$.
\end{lemma}
\begin{proof}
Suppose $m=7n$ is squarefree and is represented by $Q_2$. This forces $x \equiv y \equiv 0 \pmod{7}$. Substituting $x=7X$ and $y=7Y$ and simplifying gives $$7X^2+14Y^2+z^2 = n$$ which means that $n$ is a quadratic residue modulo $7$.
\end{proof}
We now proceed to show that if $m \not\equiv 2 \pmod{16}$ is squarefree, then $m$ is represented by $Q_1$. Again, we proceed by caess. 
\begin{itemize}
\item[(Case 1)] Suppose $m \equiv 1 \pmod{4}$. Choose a prime $a$ such that $(\tfrac{-a}{p})=1$ for all $p \vert m$, where additionally $a \equiv 5 \pmod{8}$ and $a \equiv 3 \pmod{49}$. This ensures that $$\left( \dfrac{-14m}{a} \right) = \left( \dfrac{-1}{a}\right) \left( \dfrac{2}{a}\right) \left( \dfrac{7}{a}\right) \left( \dfrac{m}{a}\right)=1.$$ Though the equation to consider now is $ab-h^2=14m$, this now behaves identically to the proof of Theorem 1 (Case 1).
\item[(Case 2)] Suppose $m \equiv 3 \pmod{4}$. Choose a prime $a$ such that $(\tfrac{-a}{p})=1$ for all $p \vert m$, where additionally $a \equiv 1 \pmod{8}$ and $a \equiv 3 \pmod{49}$. This ensures that $$\left( \dfrac{-14m}{a} \right) = \left( \dfrac{-1}{a}\right) \left( \dfrac{2}{a}\right) \left( \dfrac{7}{a}\right) \left( \dfrac{m}{a}\right)=1.$$ This now behaves like (Case 1).
\item[(Case 3)] Let $m \equiv 6, 14 \pmod{16}$. Then $m=2m'$ where $m' \equiv 3, 7,11 ,15\pmod{16}$. Let $a \equiv 1 \pmod{8}$ and $a \equiv 3 \pmod{49}$, and $(\tfrac{-a}{p})=1$ for all $p \vert m'$. This is enough to guarantee that $$\left( \dfrac{-14m}{a} \right) = \left( \dfrac{-1}{a}\right) \left( \dfrac{2}{a}\right)^2 \left( \dfrac{7}{a}\right) \left( \dfrac{m'}{a}\right)=1.$$ The proof now follows (Case 1)
\item[(Case 4)] Let $m \equiv 10 \pmod{16}$. Then $m=2m'$ where $m' \equiv 5,13 \pmod{16}$. We note that the total number of primes $p \vert m$ which are congruent to either $5$ or $7 \pmod{8}$ is odd. With that, we take $a=2a'$ where $a'$ is prime, satisfying $a' \equiv 1 \pmod{8}, a' \equiv 26 \pmod{49}$, and $\left( \tfrac{-2a'}{p} \right)=1$ for all $p \vert m'$. This gives  $$\left( \dfrac{-14m}{a'} \right) = \left( \dfrac{-1}{a'}\right) \left( \dfrac{2}{a'}\right)^2 \left( \dfrac{7}{a'}\right) \left( \dfrac{m'}{a'}\right)=1.$$ Moreover, we note that $2a' \equiv 3 \pmod{49}$, which means next considering $2a'b-h^2=14m$ we are reduced to earlier cases.
\end{itemize}

\section{Proof of Theorem 4}
We begin by noting there are five forms of determinant $23$: $Q_1: 2x^2+2xy+3y^2+2yz+5z^2$, $Q_2: x^2+y^2+23z^2$, $Q_3: x^2+2y^2+2yz+12z^2$, $Q_4: x^2+3y^2+2yz+8z^2$, and $Q_5:x^2+4y^z+2yz+6z^2$.
\begin{lemma}
If $m \equiv 1 \pmod{8}$ then $Q_1$ does not represent $m$.
\end{lemma}
\begin{proof}
Left to reader.
\end{proof}
\begin{lemma}
Let $m \in \mathbb N$ be odd. Then $Q_1$ represents $m$ if and only if $Q_1$ represents $4m$.
\end{lemma}
\begin{proof}
One direction is trivial, so suppose $Q_1$ represents $4m$ where $m$ is odd. Then $4m \equiv 4 \pmod{8}$ and a computer search will verify that in this case each of $x,y,z$ must be even in order for $Q_1(x,y,z) \equiv 4 \pmod{8}$.
\end{proof}
\begin{lemma}
If $m=23n$ where $n$ is a quadratic nonresidue modulo $23$, then $m$ is not represented by $Q_2$.
\end{lemma}
\begin{proof}
Considering $x^2+y^2+23 z^2 \equiv 0 \pmod{23}$ immediately yields $x \equiv y \equiv 0 \pmod{23}$. Substituting $x=23X$, $y=23Y$ gives 
\begin{eqnarray*}
(23X)^2 + (23Y)^2 + 23z^2 & = & 23n \\
23X^2 + 23Y^2 + z^2 & = & n,
\end{eqnarray*}
which means $n$ is a quadratic residue $\pmod{23}$.
\end{proof}
\begin{lemma}
If $m=23n$ where $n$ is a quadratic nonresidue modulo $23$, then $m$ is not represented by $Q_3$.
\end{lemma}
\begin{proof}
Setting $Q_3(x,y,z) \equiv 0 \pmod{23}$ immediately gives $x \equiv 0 \pmod{23}$. There are additional constraints on $y$ and $z$ modulo $23$. These cases behave similarly to those in previous sections, and so in the interest of space, we simply provide a summary of the data.
\begin{center}
\begin{tabular}{|c|c|}
\hline
$(y \pmod{23}, z \pmod{23})$ & $n \pmod{23}$\\
\hline
$(\pm 1, \pm 21)$ & $2$ \\
\hline
$(\pm 2, \pm 19)$ & $8$ \\
\hline
$(\pm 3, \pm 17)$ & $18$ \\
\hline
$(\pm 4, \pm 15)$ & $9$\\
\hline
$(\pm 5, \pm 13)$ & $4$\\
\hline
$(\pm 6, \pm 11)$ & $3$\\
\hline
$(\pm 7, \pm 9)$ & $6$\\
\hline
$(\pm 8, \pm 7)$ & $13$ \\
\hline
$(\pm 9, \pm 5)$ & $1$ \\
\hline
$(\pm 10, \pm 3)$ & $16$ \\
\hline
$(\pm 11, \pm 1)$ & $12$\\
\hline
\end{tabular}
\end{center}
In each case, $n$ is a quadratic residue $\pmod{23}$, which completes the proof.
\end{proof}
\begin{lemma}
If $m=23n$ where $n$ is a quadratic nonresidue modulo $23$, then $m$ is not represented by $Q_4$.
\end{lemma}
\begin{proof}
Setting $Q_4(x,y,z) \equiv 0 \pmod{23}$ immediately gives $x \equiv 0 \pmod{23}$. There are additional constraints on $y$ and $z$ modulo $23$. These cases behave similarly to those in previous sections, and so in the interest of space, we simply provide a summary of the data.
\begin{center}
\begin{tabular}{|c|c|}
\hline
$(y \pmod{23}, z \pmod{23})$ & $n \pmod{23}$\\
\hline
$(\pm 1, \pm 20)$ & $3$ \\
\hline
$(\pm 2, \pm 17)$ & $12$ \\
\hline
$(\pm 3, \pm 14)$ & $4$ \\
\hline
$(\pm 4, \pm 11)$ & $2$\\
\hline
$(\pm 5, \pm 8)$ & $6$\\
\hline
$(\pm 6, \pm 5)$ & $16$\\
\hline
$(\pm 7, \pm 2)$ & $9$\\
\hline
$(\pm 8, \pm 22)$ & $8$ \\
\hline
$(\pm 9, \pm 19)$ & $13$ \\
\hline
$(\pm 10, \pm 16)$ & $1$ \\
\hline
$(\pm 11, \pm 13)$ & $18$\\
\hline
\end{tabular}
\end{center}
In each case, $n$ is a quadratic residue $\pmod{23}$, which completes the proof.
\end{proof}
\begin{lemma}
If $m=23n$ where $n$ is a quadratic nonresidue modulo $23$, then $m$ is not represented by $Q_5$. 
\end{lemma}
\begin{proof}
Setting $Q_5(x,y,z) \equiv 0 \pmod{23}$ immediately gives $x \equiv 0 \pmod{23}$. There are additional constraints on $y$ and $z$ modulo $23$. These cases behave similarly to those in previous sections, and so in the interest of space, we simply provide a summary of the data.
\begin{center}
\begin{tabular}{|c|c|}
\hline
$(y \pmod{23}, z \pmod{23})$ & $n \pmod{23}$\\
\hline
$(\pm 1, \pm 19)$ & $4$ \\
\hline
$(\pm 2, \pm 15)$ & $16$ \\
\hline
$(\pm 3, \pm 11)$ & $13$ \\
\hline
$(\pm 4, \pm 7)$ & $18$\\
\hline
$(\pm 5, \pm 3)$ & $8$\\
\hline
$(\pm 6, \pm 22)$ & $6$\\
\hline
$(\pm 7, \pm 18)$ & $12$\\
\hline
$(\pm 8, \pm 14)$ & $3$ \\
\hline
$(\pm 9, \pm 10)$ & $2$ \\
\hline
$(\pm 10, \pm 6)$ & $9$ \\
\hline
$(\pm 11, \pm 2)$ & $1$\\
\hline
\end{tabular}
\end{center}
In each case, $n$ is a quadratic residue $\pmod{23}$, which completes the proof.
\end{proof}
Now suppose $m \not\equiv 1 \pmod{8}$ is squarefree. We will show that $m$ is represented by $Q_1$ with the following cases:
\begin{itemize}
\item[(Case 1)] Let $m \equiv 3 \pmod{4}$. Let $a$ be a prime satisfying $a \equiv 1 \pmod{4}, a \equiv 5 \pmod{529}$, and $(\frac{-a}{p})=1$ for all primes $p \vert m$. This will ensure $$\left( \dfrac{-23m}{a} \right) = \left( \dfrac{-1}{a} \right) \left( \dfrac{23}{a} \right) \left( \dfrac{m}{a} \right) = \left( \dfrac{-a}{m} \right) =1.$$ Considering the equation $ab-h^2 =23m$ and replacing $h$ with $-h$ as necessary we see that modulo $23$, $(b,h) \in \{ (0,0), (5, 5), (7, 9), (10, 2), (11, 3), (14, 1), (15, 11), (17, 4), (19, 7), (20, 10), (21,  6), (22, 8)\}.$ This automatically guarantees there is a solution to $$(A+B)^2+a+b+2h \equiv 0 \pmod{23}$$ which means that $Q$ will represent multiples of $23$. Moreover, for each case there are solutions to each of $$(A+B)^2+a+b+2h \equiv 23k$$ for $k=0,1,2,...,22$. This will guarantee that when $x \equiv y \equiv 1 \pmod{529}$ and $z=0$ that $Q$ will represent $23n$ where $n$ is a quadratic nonresidue modulo $23$ (with different congruence conditions necessary when $m$ is or is not a quadratic residue $\pmod{23}$).
\item[(Case 2)] Let $m \equiv 6 \pmod{8}$, so $m=2m'$ where $m' \equiv 3 \pmod{4}$. Let $a$ be a prime satisfying $a \equiv 1 \pmod{8}$, $a \equiv 5 \pmod{529}$, and $(\frac{-a}{p})=1$ for all primes $p \vert m'$. This yields $$\left( \dfrac{-23m}{a} \right) = \left( \dfrac{-1}{a} \right) \left( \dfrac{23}{a} \right) \left( \dfrac{2}{a} \right) \left( \dfrac{m'}{a} \right) = \left( \dfrac{-a}{m'} \right) =1.$$ And as the conditions $\pmod{23}$ on $a$ are the same as in (Case 1), the rest of the proof follows similarly.
\item[(Case 3)] Let $m \equiv 5 \pmod{8}$. We note that the total number of primes $p \vert m$ which are congruent to either $5$ or $7 \pmod{8}$ is odd. With that, we write $a=2a'$ where $a'$ is a prime satisfying $a \equiv 1 \pmod{4}$, $a \equiv 267 \pmod{529}$ and $(\frac{-a}{p})=1$ for all $p \vert m$. This first ensures $$\left( \dfrac{-23m}{a} \right) = \left( \dfrac{-1}{a} \right) \left( \dfrac{23}{a} \right) \left( \dfrac{m}{a} \right) = \left( \dfrac{-a}{m} \right) =1.$$ But also, noting that $a \equiv 2 \cdot 267 \equiv 5 \pmod{529}$ we are able to mimic previous cases at this point.
\item[(Case 4)] Let $m \equiv 2 \pmod{8}$. Let $m \equiv 2 \pmod{8}$, so $m=2m'$ where $m' \equiv 1 \pmod{4}$. Let $a$ be a prime satisfying $a \equiv 5 \pmod{8}$, $a \equiv 5 \pmod{529}$, and $(\frac{-a}{p})=1$ for all primes $p \vert m'$. This forces $$\left( \dfrac{-23m}{a} \right) = \left( \dfrac{-1}{a} \right) \left( \dfrac{23}{a} \right) \left( \dfrac{2}{a} \right) \left( \dfrac{m'}{a} \right) = \left( \dfrac{-a}{m'} \right) =1.$$ And as the conditions $\pmod{23}$ on $a$ are the same as in (Case 1), the rest of the proof follows similarly.
\end{itemize}

\section{Examples}
We end this paper with hopefully helpful if not entertaining to the reader concrete examples of choices of $A,B,a,b,h$ as outlined in the proofs of the theorems.
\begin{example}
To show that $m = 51 = 3 \cdot 17$ is represented by $2x^2+2xy+2xz+2y^2+2yz+3z^2$, we consider (Case 1) of the proof of Theorem $1$.\\
\\
We choose a prime $a \equiv 1 \pmod{4}$ and $a \equiv 3 \pmod{49}$, and (without loss of generality) $a \equiv 2 \pmod{3}$ and $ a\equiv 1 \pmod{17}$. The smallest prime satisfying all of these conditions is $a=4217$. Then solving $4217b-h^2=7 \cdot 51$ for $b$ and $h$, we see we can take $b=1613$ and $h=2608$. Note this is not the ``smallest'' solution with respect to $b>0$; however, here $b \equiv 3 \pmod{7}$ and $h \equiv 4 \pmod{7}$. Noting now that $A^2+a \equiv 0 \pmod{51}$ means modulo $51$, $A \in \{ 4,13,38,47\}$. Similarly with $B^2+b \equiv 0 \pmod{51}$ we see $B \in \{11,23,28,40 \}$. Noting, however, that we must also have $2AB+2h \equiv 0 \pmod{51}$ we see the possible pairs of $(A,B)$ modulo $51$ are $(A,B) \in \{ (4,11), (13,23), (38,28), (47,40)  \}$. Accounting for $4A^2+4AB+B^2 +4a+b+4h \equiv 2^2(A^2+a) + 2(2AB+2h)+B^2+b \equiv 21 \pmod{49}$ gives $98$ choices for $(A,B) \pmod{49}$. Among these choices is $A \equiv 0 \pmod{49}$ and $B \equiv 19 \pmod{49}$. Selecting from our $\pmod{51}$ conditions $A \equiv 38 \pmod{51}$ and $B \equiv 28 \pmod{51}$ and using the Chinese Remainder Theorem gives minimum positive values of $A=1568$ and $B=1048$. This then means $$Q(x,y,z) = 48291x^2+64544xy+3136xz+21567y^2+2096yz+51z^2.$$ 
Last, we note that  $Q(2,1,0)=7 \cdot 49117$, where as $49117 \equiv 5 \pmod{7}$ means $Q$ represents $7n$ where $n$ is a quadratic nonresidue modulo $7$. We conclude that $Q$ is equivalent to $2x^2+2xy+2xz+2y^2+2yz+3z^2$.
\end{example}

\begin{example}
To show that $m=67$ is represented by $x^2+2y^2+2yz+6z^2$ we refer to (Case 3) of the proof of Theorem $2$.\\
\\
We choose a prime $a \equiv 1 \pmod{4}$ and $a \equiv 2 \pmod{121}$ and noting $(\frac{-2}{67})=1$, we also take $a \equiv 2 \pmod{67}$. The smallest prime satisfying all of these conditions is $a=170249$. Then solving $170249b-h^2=11 \cdot 67$ for $b$ and $h$, we see we can take $b = 4413$ and $h = -27410$. The requirements of $A^2+a \equiv B^2+b \equiv 2AB+2h \equiv 0 \pmod{67}$, give $A \equiv 20 \pmod{67}$ and $B \equiv 64 \pmod{67}$, or $A \equiv 47 \pmod{67}$ and $B \equiv 3 \pmod{67}$. We choose the former. Solving $(A+B)^2+a+b+2h \equiv 22 \pmod{121}$ gives, among many options, $A \equiv 60 \pmod{121}$ and $B \equiv 0 \pmod{121}$. We then take $A=4174$ and $B=3146$. This then means $$Q(x,y,z) = 262575x^2+391164xy+8348xz+147787y^2+6292yz+67z^2.$$ We note that $Q(1,1,0)=801526 = 2 \cdot 11 \cdot 36433$, and $(\tfrac{2 \cdot 36433}{11})=-1$.
\end{example}

\begin{example}
To show that $m=26= 2 \cdot 13$ is represented by $x^2+3y^2+2yz+5z^2$, we consider (Case 4) of Theorem $3$.\\
\\
We choose a prime $a' \equiv 1 \pmod{8}$ and $a' \equiv 26 \pmod{49}$ and $a' \equiv 2 \pmod{13}$. The smallest such prime is $a' = 27809$. Set $a=2a'$. Considering next the equation $ab-h^2=14 \cdot 26$, we get $b=8440$ and $h=21666$ as a possible solution. Solving $A^2+ a \equiv 0 \pmod{26}$ and $B^2+b \equiv 0 \pmod{26}$ gives $A \in \{ 10,16\} $ and $B \in \{ 6,20\}$. Taking into account we must have $2AB+2h \equiv 0 \pmod{26}$ we see the pairs $(A,B)$ are $(10,20)$ and $(16,6)$. Because $26$ is a quadratic residue mod $7$, we next solve for $A,B \pmod{49}$ such that $(A+B)^2+a+b+2h \equiv 7 \pmod{35}$. This gives $98$ pairs. One such pair is $A \equiv 1 \pmod{49}$ and $B \equiv 43 \pmod{49}$. Using the Chinese Remainder Theorem, we take $A= 1128$ ad $B=-6$. This yields $$Q(x,y,z)= 51077x^2+1146xy+2256xz+326y^2-12yz+26.$$ Last, we note that when $x=y=1$ and $z=0$, $Q(x,y,z) = 7 \cdot 7507$, and $7507 \equiv 3 \pmod{7}$.
\end{example}


\begin{thebibliography}{fszw90}
\bibitem{BL} Blackwell, S., Durham, G., Thompson, K., and Treece, T., \textit{A generalization of a method of Mordell to ternary quadratic forms}, International Journal of Number Theory, Vol. 12 No. 8 (2016), pg 2081-2105.
\bibitem{BI} Brandt, H., Intrau, O., \textit{Tabellen reduzierter positiver tern\"arer quadratischer Formen}, Abh. S\"achs. Akad. Wiss. Math.-Nat. Kl. 45 (1958), no. 4, 261
\bibitem{Cassels1}{Cassels, J.W.S., \textit{Rational Quadratic Forms}, Dover, 2008}
\bibitem{Cassels}{Cassels, J.W.S., \textit{An Introduction to the Geometry of Numbers}, Springer-Verlag, 1959.}
\bibitem{Cox} Cox, D., \textit{Primes of the form $x^2+ny^2$}, John Wiley 
\bibitem{DHist}{Dickson, L.E., \textit{History of the Theory of Numbers, Volume III: Quadratic and Higher Forms}, Dover Publications, 2012.}
\bibitem{Dickson}{Dickson, L.E., \textit{Integers represented by positive ternary quadratic forms.}, Bull. Amer. Math. Soc. 33 (1927), 63-70.}
\bibitem{Gauss}{Gauss, C.F., \textit{Besprechung des Buchs von L.A. Seeber: Untersuchungen \"Uber die Eigenschaften der positiven tern\"aren quadratischen Forem usw. G\"ottingsche Gelehrte Anzeigen}, 1831, Juli 9. Reprinted in Werke (1876), Vol. II, 188-196.}
\bibitem{DA} Gauss, C.F., \textit{Disquisitiones Arithmeticae}, trans. A.A. Clarke, Springer New York, 1986.
\bibitem{Jones}{Jones, B., \textit{The regularity of a genus of positive ternary quadratic forms}, Trans. Amer. Math. Soc. 33 (1931), 111-124.}
\bibitem{Kaplansky} Kaplansky, I., \textit{The first nontrivial genus of positive definite ternary forms}, Mathematics of Computation, Volume 64, Number 209, January 1995, pgs. 341-345.
\bibitem{Kelley} Kelley, J., \textit{Kaplansky's ternary quadratic form}, International Journal of Mathematics and Mathematical Sciences, Volume 25, Issue 5 (2001), pg. 289-292.
\bibitem{Lam}{Lam, T.Y., \textit{Introduction to Quadratic Forms over Fields}, AMS, 2005}
\bibitem{Legendre} Legendre, A.-M., \textit{Essai sur la th\`eorie des nombres}, Paris, An VI (1797-1798)
\bibitem{Mordell}{Mordell, L.J., \textit{On the representation of a number as a sum of three squares.}, Rev. Math. Pres Appl. 3 (1958), 25-27.}
\bibitem{Nebe}{\url{http://www.math.rwth-aachen.de/~Gabriele.Nebe/LATTICES/Brandt_1.html}}
\bibitem{Ramanujan} Ramanujan, S., \textit{On the expression of a number in the form $ax^2+by^2+cz^2+du^2$}, Proc. Camb. Phil. Soc. 19 (1916), 11-21.
\bibitem{Rouse} Rouse, J., \textit{Quadratic Forms Representing All Odd Positive Integers}, American Journal of Mathematics. Vol 136 (2011).
\bibitem{Sage} SageMath, the Sage Mathematics Software System (Version 9.3), The Sage Developers, 2015, http://www.sagemath.org
\end{thebibliography}
\end{document}